\documentclass[times]{article}

\usepackage{amsmath}
\usepackage{url}
\usepackage{citesort}
\usepackage[amssymb]{SIunits}
\usepackage{lscape}
\usepackage{graphicx}

\def\Xint#1{\mathchoice
   {\XXint\displaystyle\textstyle{#1}}%
   {\XXint\textstyle\scriptstyle{#1}}%
   {\XXint\scriptstyle\scriptscriptstyle{#1}}%
   {\XXint\scriptscriptstyle\scriptscriptstyle{#1}}%
   \!\int}
\def\XXint#1#2#3{{\setbox0=\hbox{$#1{#2#3}{\int}$}
     \vcenter{\hbox{$#2#3$}}\kern-.5\wd0}}
\def\ddashint{\Xint=}
\def\dashint{\Xint-}

\newcommand{\D}{\mathrm{d}}  
\newcommand{\rbar}{\overline{r}}
\newcommand{\Rbar}{\overline{R}}

\begin{document}

\bibliographystyle{unsrt}



\title{Potential integrals on triangles}

\author{Michael Carley}

\maketitle

\begin{abstract}
  The problem of evaluating potential integrals on planar triangular
  elements has been addressed using a polar coordinate
  decomposition. The resulting formulae are general, exact, easily
  implemented, and have only one special case, that of a field point
  lying in the plane of the element. Results are presented for the
  evaluation of the potential and its gradients, where the integrals
  must be treated as principal values or finite parts, for elements
  with constant and linearly varying source terms. These results are
  tested by application to a single triangular element to the
  evaluation of the potential gradient outside the unit cube. In both
  cases, the method is shown to be accurate and convergent.
\end{abstract}

\section{INTRODUCTION}
\label{sec:intro}
\vspace{-2pt}

A basic operation in any code for the Boundary Element Method (BEM) is
the evaluation of potential integrals on elements, whether in the
solution of the integral equation, or in the evaluation of field
quantities. If we consider the Laplace problem, $\phi$ the potential
external to a surface $S$ is given by the integral formulation:
\begin{align}
  \label{equ:potential}
  \phi(\mathbf{x}) &= 
  \int_{S} 
  \frac{\partial\phi_{1}}{\partial n}G(\mathbf{x},\mathbf{x}_{1})
  -
  \frac{\partial G(\mathbf{x},\mathbf{x}_{1})}{\partial n}\phi_{1}
  \,\D S,
\end{align}
where $\mathbf{x}$ indicates position, subscript $1$ variables of
integration on the surface $S$ and $n$ the outward pointing normal to
the surface. The Green's function $G$ is given by
\begin{align}
  \label{equ:laplace}
  G(\mathbf{x};\,\mathbf{x}_{1}) &= \frac{1}{4\pi R},\\
  R &= |\mathbf{x}-\mathbf{x}_{1}|. \nonumber
\end{align}
Given the surface potential $\phi$ and gradient $\partial\phi/\partial
n$, the potential, and, after differentiation, the gradient(s), can be
evaluated at any point in the field. Also, given a boundary condition
for $\phi$ and/or $\partial\phi/\partial n$ on $S$, the integral
equation can be solved for $\mathbf{x}\in S$.

In any case, the method of solution remains the same: the surface $S$
is divided into elements and suitable shape functions are used to
interpolate the potential on these elements. The integral equation is
transformed to a linear system in the element potentials, with the
influence coefficients determined by the potential generated by each
element at each node of the surface mesh. This leads to the
requirement to evaluate integrals $I$ and $\partial I/\partial n$
where: 
\begin{align}
  \label{equ:element}
  I &= \iint f(\xi,\eta)
  G(\mathbf{x},\mathbf{x}_{1}(\xi,\eta))\,\D \xi\,\D \eta,
\end{align}
where $E$ is the surface of an element and $(\xi,\eta)$ is a
coordinate system local to $E$.

Numerous techniques have been proposed for the evaluation of $I$ and
its derivatives. Many of these have been numerical and have often
focussed on the question of evaluating the integral when the
evaluation point is close to, but not on, the element~\cite[for
example]{telles87,hayami-matsumoto94,hayami05}. For the Laplace
equation, however, it is possible to exactly evaluate the potential
integrals analytically and a number of techniques have been published
using this approach~\cite[for
example]{okon-harrington82a,okon-harrington82b,okon85,newman86,%
  nintcheu-fata09,carini-salvadori02,salvadori10}, which has the
advantages of being exact and of being less prone to numerical errors
introduced by singularities and near-singularities, and, often, more
efficient than numerical quadrature.

This paper introduces a method for the evaluation of potential
integrals on flat triangular elements, in practice the most common
problem encountered in BEM integration. The motivation for this work
is to simplify the implementation of the analytical formulae. The
numerous results which have been published hereto are, obviously,
algebraically equivalent in that they give, or should give, the same
answer for any given combination of element and field point. In
practice, however, they are not numerically equivalent and, in
particular, they have different special cases which must be
handled. For example, the formulae of Newman~\cite{newman86}, based on
the use of Green's theorem to reduce the surface integral to a
sequence of line integrals on the element boundary, have a special
case when the field point is collinear with an edge of the
triangle. In Salvadori's work, there are a number of special cases
which must be identified and handled separately depending on
particular combinations of parameters~\cite[Figure~5]{salvadori10}. 

In this paper, we present formulae based on a polar coordinate
decomposition of the basic integrals which simplifies the integrands
to the point where they can be easily evaluated using standard
relations and tables and which require handling of only one special
case, that when the field point lies in the plane of the element. The
method can be generalized to source terms of any order, and to
gradients of the potential, which give rise to strongly singular
integrands which require treatment as hypersingular integrals.

\section{Potential integrals}
\label{sec:integrals}

The approach adopted in developing formulae for the evaluation of
potential integrals is similar to that used in previous work: the
integral is evaluated on a reference triangle into which the real
element can be transformed. In this case, the reference triangle lies
in the plane $z=0$, requiring a rotation of the coordinate system, and
has one vertex at the projection of the field point into that plane,
requiring a decomposition of the original element into a set of
subtriangles. We begin by developing formulae for the integrals over
the reference triangle, before explaining how these formulae can be
applied to general elements. 

\begin{figure}
  \centering
  \includegraphics{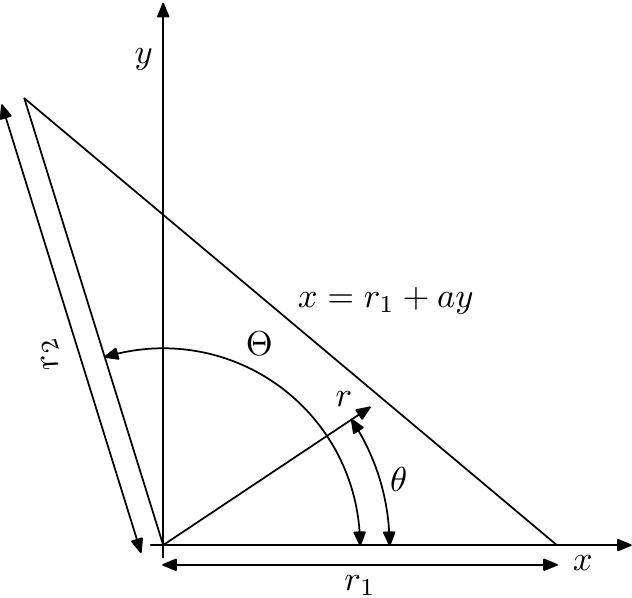}
  \caption{Reference triangle for quadrature}
  \label{fig:triangle}
\end{figure}

Figure~\ref{fig:triangle} shows the reference triangle and associated
notation. The triangle lies in the plane $z=0$, with one vertex at the
origin. The two sides which meet at the origin have length $r_{1}$ and
$r_{2}$ and subtend an angle $\Theta$. The third side is given by
$x=r_{1}+ay$ with:
\begin{align}
  \label{equ:slope}
  a &= \frac{r_{2}\cos\Theta - r_{1}}{r_{2}\sin\Theta}.
\end{align}
Integrations will be performed by transforming the double integrals in
into integrals in the polar coordinate system $(r,\theta)$ centred on
the origin. In the triangle coordinate system, the field point is
placed at $(0,0,z)$.

The integrals which are required are of the following forms, with the
factor $1/4\pi$ suppressed:
\begin{subequations}
  \label{equ:basic}
  \begin{align}
    \label{equ:gfunc}
    I &= \iint f(x_{1},y_{1})\frac{1}{R}\,\D x_{1}\,\D y_{1},\\
    \label{equ:gfunc:z}
    \frac{\partial I}{\partial z}
    &=
    -z \iint  f(x_{1},y_{1})\frac{1}{R^{3}}\,\D x_{1}\,\D y_{1},\\
    \label{equ:gfunc:x}
    \frac{\partial I}{\partial x}
    &=
    -\iint  f(x_{1},y_{1})\frac{x-x_{1}}{R^{3}}\,\D x_{1}\,\D y_{1},\\
    \label{equ:gfunc:zz}
    \frac{\partial^{2} I}{\partial z^{2}}
    &=
    \iint  
    f(x_{1},y_{1})
    \left[
      3\frac{z^{2}}{R^{5}} -\frac{1}{R^{3}}
    \right]
    \,\D x_{1}\,\D y_{1},\\
    \label{equ:gfunc:xz}
    \frac{\partial^{2} I}{\partial x\partial z}
    &=
    \iint  
    f(x_{1},y_{1})
    \frac{3z(x-x_{1})}{R^{5}}
    \,\D x_{1}\,\D y_{1},\\
    \label{equ:gfunc:xx}
    \frac{\partial^{2} I}{\partial x^{2}}
    &=
    \iint  
    f(x_{1},y_{1})
    \left[
      3\frac{(x-x_{1})^{2}}{R^{5}} -\frac{1}{R^{3}}
    \right]
    \,\D x_{1}\,\D y_{1},\\
    R^{2} &= (x-x_{1})^{2} + (y-y_{1})^{2} + z^{2}.\nonumber
  \end{align}
\end{subequations}
In these integrals, $f(x_{1},y_{1})$ is a source term which varies
over the element and subscript $1$ denotes a variable of
integration. The basic integral $I$ is the integral over a general
triangle in the plane $z=0$ of the Green's function weighted on the
source term. Differentiation with respect to $z$,
Equation~\ref{equ:gfunc:z}, corresponds to taking the normal
derivative and derivatives with respect to $(x,y,z)$ are used to find
the gradient of potential, e.g. the velocity in potential flow
problems. Further differentiation yields the hypersingular formulation
used in certain applications to avoid numerical difficulties
introduced by spurious eigenvalues~\cite{burton-miller71}.

The integrals of Equations~\ref{equ:gfunc} can be evaluated on the
reference triangle of Figure~\ref{fig:triangle} as
follows. Considering a general integral:
\begin{align}
  \label{equ:basic:1}
  I &= \iint f(x_{1},y_{1})R^{\gamma}\,\D x_{1}\,\D y_{1},\\
  \label{equ:basic:2}
  &= \int_{0}^{\Theta} \int_{0}^{\rbar(\theta)}
  f(x_{1},y_{1})R^{\gamma}r\,\D r\D \theta,
  \quad
  R^{2} = r^{2}+z^{2},
\end{align}
where
\begin{align}
  \label{equ:rbar}
  \rbar(\theta) &=
  \frac{r_{1}}{(1+a^{2})^{1/2}}\frac{1}{\cos(\theta+\phi)},\quad
  \phi = \tan^{-1} a,
\end{align}
$I$ can be evaluated analytically using tabulated relations. 

The first case considered is that of an element with constant source
distribution $f(x_{1},y_{1})\equiv1$:
\begin{align}
  \label{equ:constant:1}
  I &= \int_{0}^{\Theta} \int_{0}^{\rbar(\theta)}R^{\gamma}r\,\D r\D \theta.
\end{align}
Under a change of variables:
\begin{align}
  \label{equ:constant:2}
  I &= \int_{0}^{\Theta} \int_{|z|}^{\Rbar(\theta)}R^{\gamma+1} \D R\D
  \theta,  \\
  &= \frac{1}{\gamma+2} \int_{0}^{\Theta} \Rbar^{\gamma+2} -
  |z|^{\gamma+2}\,\D\theta,\\
  \Rbar^{2} &= \rbar^{2} + z^{2},\nonumber
\end{align}
which can be evaluated with a further change of variables from
$\theta$ to $\theta+\phi$:
\begin{align}
  \label{equ:constant:3}
  I &= \frac{1}{\gamma+2}
  \left[
    \beta^{\gamma+2}
    \int_{\phi}^{\Theta+\phi}
    \frac{\Delta^{\gamma+2}}{\cos^{\gamma+2}\theta}
    \,\D\theta
    -
    |z|^{\gamma+2}\Theta
  \right],
\end{align}
where:
\begin{align*}
  \Delta^{2} &= 1 - \alpha^{2}\sin^{2}\theta,\quad
  \alpha^{2} = z^{2}/\beta^{2},\quad
  \alpha' = \left(1-\alpha^{2}\right)^{1/2},\quad
  \beta^{2} = \frac{r_{1}^{2}+z^{2}(1+a^{2})}{1+a^{2}},\nonumber\\
  \rbar &= \frac{\beta\alpha'}{\cos\theta},\quad\text{and}\quad
  \Rbar = \beta\frac{\Delta}{\cos\theta}.
\end{align*}
Integrals of the form of Equation~\ref{equ:constant:3} are well
understood and extensively
tabulated~\cite{gradshteyn-ryzhik80,prudnikov-brychkov-marichev03:1}. To
simplify the statement of future results, we write the integral as:
\begin{align}
  \label{equ:constant:4}
  I &= \frac{1}{\gamma+2}
  \left[
    \beta^{\gamma+2}
    \left(
      I_{0,-\gamma-2}^{(\gamma+2)}(\alpha,\Theta+\phi) -
      I_{0,-\gamma-2}^{(\gamma+2)}(\alpha,\phi)
    \right)
    -|z|^{\gamma+2}\Theta
  \right],
\end{align}
where
\begin{align}
  \label{equ:Imn}
  I_{m,n}^{(p)}(\alpha,\theta) &= 
  \int_{0}^{\theta}\sin^{m}\theta'\cos^{n}\theta' \Delta^{p}\,\D\theta'.
\end{align}
This result is the general form for the integral of any integer power
of $R$ over the reference triangle. In particular, subject to triangle
decomposition and rotation, it is the only integral required in a BEM
which employs constant source elements, as it yields the integrals of
$G$ and $\partial G/\partial n$, and includes the integral of
$1/R^{3}$ considered by~\cite{carini-salvadori02}. The only special
case which need be handled is that for $z=0$, when the field point
lies in the element plane, which will be discussed in
Section~\ref{sec:in:plane}.

Higher order elements can be handled in a similar fashion. In the case
of elements with linear variation in $f(x_{1},y_{1})$, the integrals
to be evaluated are of the form:
\begin{align}
  \label{equ:linear}
  I &= \int_{\phi}^{\Theta+\phi} 
  \int_{0}^{\rbar(\theta)}
  \frac{r\sin^{m}\theta\cos^{n}\theta}{R}r\,\D r\,\D\theta,\quad
  m+n=1.
\end{align}
In the case $m=0,n=1$, this yields:
\begin{align}
  \label{equ:linear:1}
  I &= \int_{\phi}^{\Theta+\phi} \cos\theta 
  \int_{0}^{\rbar(\theta)} \frac{r^{2}}{R}\,\D r\,\D\theta,\\
  \label{equ:linear:2}
  &=
  \frac{\beta^{2}\alpha'}{2}
  \int_{\phi}^{\Theta+\phi}\frac{\Delta}{\cos\theta}\,\D\theta
  -
  \frac{z^{2}}{4}
  \left.
    \log\frac{\Delta+\alpha'}{\Delta-\alpha'}\sin\theta
  \right|_{\phi}^{\Theta+\phi}
  +
  \frac{z^{2}\alpha'}{2}
  \int_{\phi}^{\Theta+\phi}\frac{\sin^{2}\theta}{\cos\theta}
  \frac{1}{\Delta}\,\D\theta,\\
  &= 
  \frac{\beta^{2}\alpha'}{2}
  \left[
    I_{0,-1}^{(1)}(\alpha,\Theta+\phi) - 
    I_{0,-1}^{(1)}(\alpha,\phi)
  \right]
  \nonumber\\
  &+
  \label{equ:linear:2a}
  \frac{z^{2}\alpha'}{2}
  \left[
    I_{2,-1}^{(-1)}(\alpha,\Theta+\phi) - 
    I_{2,-1}^{(-1)}(\alpha,\phi)
  \right]
  -\frac{z^{2}}{4}
  \left.
    \log\frac{\Delta+\alpha'}{\Delta-\alpha'}\sin\theta
  \right|_{\phi}^{\Theta+\phi},
\end{align}
upon integration by parts. Likewise,
\begin{align}
  \int_{\phi}^{\Theta+\phi} \sin\theta
  \int_{0}^{\rbar(\theta)} \frac{r^{2}}{R}\,\D r\,\D\theta
  &=
  \frac{\beta^{2}\alpha'}{2}
  \left[
    I_{1,-2}^{(1)}(\alpha,\Theta+\phi) - 
    I_{1,-2}^{(1)}(\alpha,\phi)
  \right]
  \nonumber\\
  \label{equ:linear:3}
  &-
  \frac{z^{2}\alpha'}{2}
  \left[
    I_{1,0}^{(-1)}(\alpha,\Theta+\phi) - 
    I_{1,0}^{(-1)}(\alpha,\phi)
  \right]
  +\frac{z^{2}}{4}
  \left.
    \log\frac{\Delta+\alpha'}{\Delta-\alpha'}\cos\theta
  \right|_{\phi}^{\Theta+\phi}.
\end{align}

Finally, higher order elements can be reduced to combinations of
the integrals for constant and linear source triangles using the
identity $r^{2}=R^{2}-z^{2}$. For second order elements:
\begin{align}
  I &=
  \int_{\phi}^{\Theta+\phi} \sin^{m}\theta\cos^{n}\theta
  \int_{0}^{\rbar}R^{\gamma}r^{3}\,\D r\,\D \theta, \quad m+n=2.
  \nonumber \\
  &=
  \label{equ:second:1}
  \int_{\phi}^{\Theta+\phi} \sin^{m}\theta\cos^{n}\theta
  \int_{0}^{\rbar}R^{\gamma}(R^{2}-z^{2})r\,\D r\,\D \theta,
\end{align}
which reduces to the standard integrals previously introduced:
\begin{align}
  \int_{\phi}^{\Theta+\phi} \sin^{m}\theta\cos^{n}\theta
  \int_{0}^{\rbar}R^{\gamma}r^{3}\,\D r\,\D \theta &=
  \frac{\beta^{\gamma+4}}{\gamma+4}
  \left[
    I_{m,n-\gamma-4}^{(\gamma+4)}(\alpha,\Theta+\phi) -
    I_{m,n-\gamma-4}^{(\gamma+4)}(\alpha,\phi)     
  \right]\nonumber\\
  &-
  z^{2}\frac{\beta^{\gamma+2}}{\gamma+2}
  \left[
    I_{m,n-\gamma-2}^{(\gamma+2)}(\alpha,\Theta+\phi) -
    I_{m,n-\gamma-2}^{(\gamma+2)}(\alpha,\phi)     
  \right]\nonumber \\
  \label{equ:second:2}
  &+
  |z|^{\gamma+4}
  \left(
    \frac{1}{\gamma+2}
    -
    \frac{1}{\gamma+4}
  \right)
  \left[J_{m,n}(\Theta+\phi)-J_{m,n}(\phi)\right],\\
  \label{equ:jfunc}
  J_{m,n}(\theta) &=
  \int_{0}^{\theta}\sin^{m}\theta'\cos^{n}\theta'\,\D\theta'.
\end{align}
Explicit expressions for the main results required for constant and
linear elements are given in Table~\ref{tab:results}, with
$I_{m,n}^{(p)}$ and $J_{m,n}$ given in the appendix,
Tables~\ref{tab:Imnp} and~\ref{tab:Jmn}.

\subsection{Field point in plane}
\label{sec:in:plane}

The formulae presented above are valid for all field points out of the
triangle plane. The case of a point lying in the plane---the only
special case which arises---must be handled separately for two
reasons. The first is that the formulae presented so far break down
numerically unless certain limits are taken explicitly; the second is
that a number of the integrals in question are truly singular and must
be interpreted as principal value or hypersingular integrals. It is
more convenient to handle all of the in-plane cases together using the
approach of Brand\~{a}o~\cite{brandao87} who gives a convenient
analysis for singular integrals of general form. For example, taking
the integral of $1/R^{3}$, and noting that for $z=0$, $r=R$:
\begin{align}
  \label{equ:in:plane:1}
  \int_{0}^{\Theta}
  \int_{0}^{\rbar}
  \frac{1}{R^{3}}
  r\,\D r\,
  \D\theta &=
  \int_{0}^{\Theta}
  \int_{0}^{\rbar}
  \frac{1}{r^{2}}
  \,\D r\,
  \D\theta.
\end{align}
Interpreting the integral in $r$ as a finite-part integral, using
Brand\~{a}o's method, it can be written:
\begin{align}
  \label{equ:in:plane:2}
  \ddashint_{0}^{\rbar}
  \frac{1}{r^{2}}
  \D r
  &=
  -\frac{1}{\rbar} = -\frac{\cos\theta}{\beta},
\end{align}
so that:
\begin{align}
  \label{equ:in:plane:3}
  \int_{0}^{\Theta}
  \ddashint_{0}^{\rbar}
  \frac{1}{R^{3}}
  \D r\,\D\theta
  &=
  -\frac{1}{\beta}
  \left[
    J_{0,1}(\Theta+\phi) - J_{0,1}(\phi)
  \right].
\end{align}
For a linear source term:
\begin{align}
  \label{equ:in:plane:4}
  \int_{\phi}^{\Theta+\phi}
  \int_{0}^{\rbar}
  \frac{r\cos\theta}{R^{3}}r
  \,\D r\,\D\theta
  &=
  \int_{\phi}^{\Theta+\phi}
  \int_{0}^{\rbar}
  \frac{\cos\theta}{r}\,
  \D r\,\D\theta.
\end{align}
The result of Brand\~{a}o's which is applicable is:
\begin{align}
  \label{equ:in:plane:5}
  \dashint_{0}^{\rbar}\frac{1}{r}\,\D r
  &=
  \log\rbar = -\log\frac{\cos\theta}{\beta},
\end{align}
so that, upon integration by parts:
\begin{align}
  \label{equ:in:plane:6}
  \int_{0}^{\Theta}
  \dashint_{0}^{\rbar}
  \frac{r\cos(\theta+\phi)}{R^{3}}\,
  \D r\,\D\theta
  &=
  \left.
    -\sin\theta\log\frac{\cos\theta}{\beta}
  \right|_{\phi}^{\Theta+\phi} - 
  \left[
    J_{2,-1}(\Theta+\phi) - J_{2,-1}(\phi)
  \right]
\end{align}
and, similarly:
\begin{align}
  \label{equ:in:plane:7}
  \int_{0}^{\Theta}
  \dashint_{0}^{\rbar}
  \frac{r\sin(\theta+\phi)}{R^{3}}\,
  \D r\,\D\theta
  &=
  \left.
    \cos\theta\log\frac{\cos\theta}{\beta}
  \right|_{\phi}^{\Theta+\phi} +
  \left[
    J_{1,0}(\Theta+\phi) - J_{2,-1}(\phi)
  \right]
\end{align}
Finally, for integrals containing terms $R^{-5}$, the corresponding
results are:
\begin{align}
  \label{equ:hyper:4}
  \int_{0}^{\Theta}
  \ddashint_{0}^{\rbar}
  \frac{1}{R^{5}}
  \D r\,\D\theta
  &=
  -\frac{1}{3\beta^{3}}
  \left[
    J_{0,3}(\Theta) - J_{0,3}(\phi)
  \right]
  \\
  \label{equ:hyper:5}
  \int_{0}^{\Theta}
  \ddashint_{0}^{\rbar}
  \frac{r\cos(\theta+\phi)}{R^{5}}
  \D r\,\D\theta
  &=
  -\frac{1}{2\beta^{2}}
  \left[
    J_{0,3}(\Theta) - J_{0,3}(\phi)
  \right]  
  \\
  \label{equ:hyper:6}  
  \int_{0}^{\Theta}
  \ddashint_{0}^{\rbar}
  \frac{r\sin(\theta+\phi)}{R^{5}}
   \D r\,\D\theta
  &=
  -\frac{1}{2\beta^{2}}
  \left[
    J_{1,2}(\Theta) - J_{1,2}(\phi)
  \right]  
\end{align}

\begin{table}
  \caption{Summary of main results for constant and linear triangular
    elements}
  \label{tab:results}
  \centering
  
  \begin{tabular}{crr}
    \hline
    \multicolumn{1}{c}{Integrand} &
    \multicolumn{1}{c}{Integral, $z\neq0$} &
    \multicolumn{1}{c}{Integral, $z=0$} \\
    \hline
    $\displaystyle\frac{1}{R}$ 
    & 
    $\displaystyle \beta I_{0,-1}^{(1)} - |z|J_{0,0}$
    &
    $\displaystyle \beta J_{0,-1}$
    \\[2ex]
    $\displaystyle\frac{r\cos(\theta+\phi)}{R}$ 
    & 
    $\displaystyle \frac{\beta^{2}\alpha'}{2} I_{0,-1}^{(1)} +
    \frac{z^{2}\alpha'}{2} I_{2,-1}^{(-1)}
    - \frac{z^{2}}{4}\sin\theta\log\frac{\Delta+\alpha'}{\Delta-\alpha'}$
    &
    $\displaystyle \frac{\beta^{2}}{2}J_{0,-1}$
    \\[2ex]
    $\displaystyle\frac{r\sin(\theta+\phi)}{R}$ 
    & 
    $\displaystyle \frac{\beta^{2}\alpha'}{2}I_{1,-2}^{(1)} - 
    \frac{z^{2}\alpha'}{2} I_{1,0}^{(-1)}
    + \frac{z^{2}}{4}\cos\theta\log\frac{\Delta+\alpha'}{\Delta-\alpha'}$
    &
    $\displaystyle \frac{\beta^{2}}{2}J_{1,-2}$
    \\[2ex]
    $\displaystyle \frac{1}{R^{3}}$ 
    & 
    $\displaystyle -\frac{1}{\beta}I_{0,1}^{(-1)} +
    \frac{1}{|z|}J_{0,0}$
    &
    $\displaystyle -\frac{1}{\beta}J_{0,1}$
    \\[2ex]
    $\displaystyle\frac{r\cos(\theta+\phi)}{R^{3}}$ 
    & 
    $\displaystyle -\alpha'I_{0,1}^{(-1)} -
    \alpha' I_{2,-1}^{(-1)} + 
    \frac{\sin\theta}{2}\log\frac{\Delta+\alpha'}{\Delta-\alpha'}
    $
    &
    $\displaystyle -\sin\theta\log\frac{\cos\theta}{\beta} - J_{2,-1}$
    \\[2ex]
    $\displaystyle\frac{r\sin(\theta+\phi)}{R^{3}}$ 
    & 
    $\displaystyle
    -\frac{\cos\theta}{2}\log\frac{\Delta+\alpha'}{\Delta-\alpha'}
    $
    &
    $\displaystyle \cos\theta\log\frac{\cos\theta}{\beta} + J_{1,0}$
    \\[2ex]
    $\displaystyle \frac{1}{R^{5}}$ 
    & 
    $\displaystyle 
    -\frac{1}{3}
    \left[
      \frac{1}{\beta^{3}}I_{0,3}^{(-3)} - \frac{1}{|z|^{3}}J_{0,0}
    \right]
    $
    &
    $\displaystyle -\frac{1}{3\beta^{3}}J_{0,3}$
    \\[2ex]
    $\displaystyle \frac{r\cos(\theta+\phi)}{R^{5}}$ 
    & 
    $\displaystyle 
    \frac{\alpha'^{3}}{3z^{2}}I_{0,1}^{(-3)}
    $
    &
    $\displaystyle -\frac{1}{2\beta^{2}}J_{0,3}$
    \\[2ex]
    $\displaystyle \frac{r\sin(\theta+\phi)}{R^{5}}$ 
    & 
    $\displaystyle 
    \frac{\alpha'^{3}}{3z^{2}}I_{1,0}^{(-3)}
    $
    &
    $\displaystyle -\frac{1}{2\beta^{2}}J_{1,2}$
    \\[2ex]
    $\displaystyle \frac{r^{2}\cos^{2}(\theta+\phi)}{R^{3}}$ 
    & 
    $\displaystyle 
    \beta I_{0,1}^{(1)} + \frac{z^{2}}{\beta}I_{0,3}^{(-1)} - 2|z|J_{0,2}
    $
    &
    $\beta J_{0,1}$
    \\[2ex]
    $\displaystyle \frac{r^{2}\cos(\theta+\phi)\sin(\theta+\phi)}{R^{3}}$ 
    & 
    $\displaystyle 
    \beta I_{1,0}^{(1)} + \frac{z^{2}}{\beta}I_{1,2}^{(-1)} - 2|z|J_{1,1}
    $
    &
    $\beta J_{1,0}$
    \\[2ex]
    $\displaystyle \frac{r^{2}\sin^{2}(\theta+\phi)}{R^{3}}$ 
    & 
    $\displaystyle 
    \beta I_{2,-1}^{(1)} + \frac{z^{2}}{\beta}I_{2,1}^{(-1)} - 2|z|J_{2,0}
    $
    &
    $\beta J_{2,-1}$
    \\[2ex]
    $\displaystyle \frac{r^{2}\cos^{2}(\theta+\phi)}{R^{5}}$ 
    & 
    $\displaystyle 
    -\frac{1}{\beta} I_{0,3}^{(-1)} + \frac{z^{2}}{3\beta^{3}}I_{0,5}^{(-3)} 
    - \frac{2}{3|z|}J_{0,2}
    $
    &
    $\displaystyle -\frac{1}{\beta}J_{0,3}$
    \\[2ex]
    $\displaystyle \frac{r^{2}\cos(\theta+\phi)\sin(\theta+\phi)}{R^{5}}$ 
    & 
    $\displaystyle 
    -\frac{1}{\beta} I_{1,2}^{(-1)} + \frac{z^{2}}{3\beta^{3}}I_{1,4}^{(-3)} 
    - \frac{2}{3|z|}J_{1,1}
    $
    &
    $\displaystyle -\frac{1}{\beta}J_{1,2}$
    \\[2ex]
    $\displaystyle \frac{r^{2}\sin^{2}(\theta+\phi)}{R^{5}}$ 
    & 
    $\displaystyle 
    -\frac{1}{\beta} I_{2,1}^{(-1)} + \frac{z^{2}}{3\beta^{3}}I_{2,3}^{(-3)} 
    - \frac{2}{3|z|}J_{2,0}
    $
    &
    $\displaystyle -\frac{1}{\beta}J_{2,1}$
    \\
    \hline
  \end{tabular}
\end{table}

\subsection{Integration over general triangles}
\label{sec:general}

The preceding sections give exact results for integration over the
reference triangle of Figure~\ref{fig:triangle}. Any triangular
element can be transformed to a combination of triangles of reference
shape by the following procedure. It is assumed that the coordinate
system has already been transformed so that the triangle lies in the
plane $z=0$, a routine operation in evaluation of potential integrals.

\begin{figure}
  \centering
  \begin{tabular}{cc}
    \includegraphics{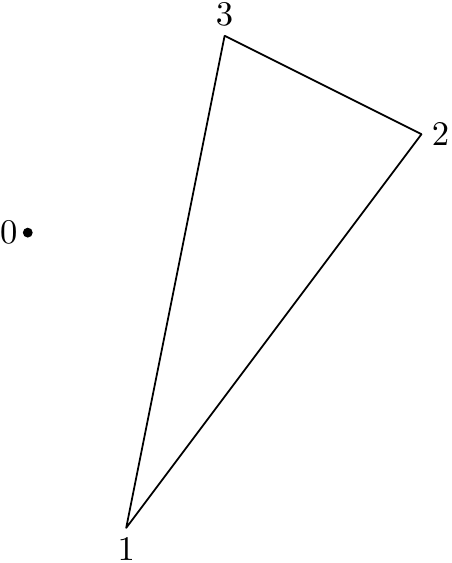} &
    \includegraphics{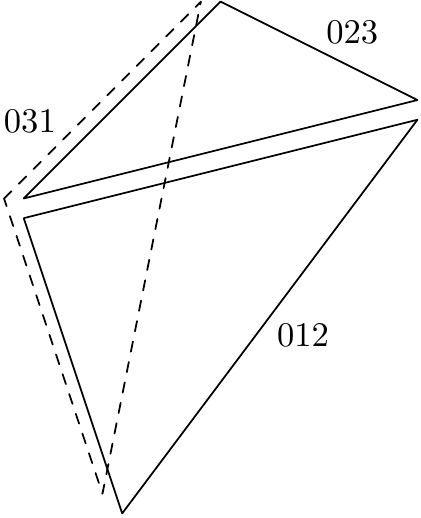} 
  \end{tabular}
  \caption{Integration over a general triangle (left) by subdivision
    into three triangles centred at the origin (right). The triangle
    shown dashed on the right has negative orientation and its
    contribution is subtracted from that of the other two.}
  \label{fig:general}
\end{figure}

Figure~\ref{fig:general} shows a general triangle $123$ and the
projection of a field point into the plane of the triangle at point
$0$. The system is decomposed into three subtriangles $012$, $023$,
and $031$ which are all of reference form. The integral over the
triangle is then the sum of the integrals over the subtriangles. In
summing the subtriangle integrals, account must be taken of the
orientation of the triangles. For example, in
Figure~\ref{fig:general}, the contributions of triangles $012$ and
$023$ are added, while that of $031$ is subtracted. The simplest way
to incorporate the orientations is to apply it to the calculation of
the angle $\Theta$ in Figure~\ref{fig:triangle}:
\begin{align}
  \label{equ:theta}
  \Theta &= 
  \pm \cos^{-1}
  \frac{(\mathbf{x}_{2}-\mathbf{x}_{0}).(\mathbf{x}_{1}-\mathbf{x}_{0})}
  {|\mathbf{x}_{2}-\mathbf{x}_{0}||\mathbf{x}_{1}-\mathbf{x}_{0}|},
\end{align}
where the sign is chosen to agree with the orientation (signed area)
of the triangle. Orientations are computed using Shewchuk's robust
adaptive predicates~\cite{shewchuk96a}.

For constant source elements, the subtriangle integrals can be summed
directly. When there is a variation in the source term, two rotations
must be applied, one by angle $\phi$, Equation~\ref{equ:rbar}, and one
by angle $\psi$:
\begin{align}
  \label{equ:psi}
  \psi = \tan^{-1}\frac{y_{1}-y_{0}}{x_{1}-x_{0}},
\end{align}
where subscript $0$ refers to the field point and subscript $1$ to the
first vertex of the triangle. Assembling these rotations gives, for
the linear source:
\begin{align}
  \left[
    \begin{array}{l}
      \displaystyle
      \iint \frac{x_{1}}{R} \D x_{1}\,\D y_{1} \\
      \displaystyle
      \iint \frac{y_{1}}{R} \D x_{1}\,\D y_{1}
    \end{array}
  \right]
  &=
  \sum_{i}
  \left[
    \begin{array}{rr}
      \cos\psi_{i} & -\sin\psi_{i} \\
      \sin\psi_{i} & \cos\psi_{i}
    \end{array}
  \right]
  \left[
    \begin{array}{rr}
      \cos\phi_{i} & \sin\phi_{i} \\
      -\sin\phi_{i} & \cos\phi_{i}
    \end{array}
  \right]
  \left[
    \begin{array}{l}
      \displaystyle
      \int_{\phi}^{\phi+\Theta}
      \int_{0}^{\rbar} 
      \frac{r\cos(\theta+\phi)}{R} r\,\D r\,\D\theta\\
      \displaystyle
      \int_{\phi}^{\phi+\Theta}
      \int_{0}^{\rbar} \frac{r\sin(\theta+\phi)}{R} r\,\D r\,\D\theta
    \end{array}
  \right] \nonumber \\
  &+
  \left[
    \begin{array}{r}
      x_{0} \\
      y_{0}
    \end{array}
  \right]
  \int_{\phi}^{\phi+\Theta}
  \int_{0}^{\rbar} \frac{1}{R} r\,\D r\,\D\theta,
  \label{equ:composition}  
\end{align}
where subscript $i=1,2,3$ refers to the subtriangles into which the
element is decomposed. 

\subsection{Summary of algorithm}
\label{sec:summary}

In summary, given a triangle lying in the plane $z=0$ and a field
point at some general position, the potential integrals can be
computed as follows:
\begin{enumerate}
\item check if the field point lies in the plane $z=0$;
\item decompose the triangle into subtriangles as in
  Figure~\ref{fig:general};
\item for each triangle $i=1,2,3$:
  \begin{enumerate}
  \item compute the orientation (signed area) of the triangle;
  \item if the orientation is not zero, calculate $a$, $r_{1}$,
    $r_{2}$, $\Theta_{i}$, and $\phi_{i}$ using
    Equations~\ref{equ:slope} and~\ref{equ:rbar}, taking account of
    the orientation;
  \item evaluate the required integrals using the appropriate formulae
    of Table~\ref{tab:results}, depending on whether or not $z=0$;
  \item if required, apply the rotations of
    Equation~\ref{equ:composition}. 
  \end{enumerate}
\item sum the results from the three subtriangles.
\end{enumerate}
If necessary, e.g. for the computation of gradients, a further
rotation can then be applied to return to the global coordinate
system.

\section{Numerical tests}
\label{sec:numerical}

Two main sets of tests have been performed to check the formulae
presented. The first consists of the evaluation of sample integrals on
a single triangle, with field points chosen to be representative of
the cases likely to cause numerical difficulties. The second is a test
of the estimation of the gradient of a known potential field outside
the unit cube. This tests the method in a BEM code and assesses its
ability to accurately integrate the strongly singular integrands which
arise in computing gradients. 

\subsection{Integrals on a sample triangle}
\label{sec:triangle}

\begin{figure}
\centering
\includegraphics{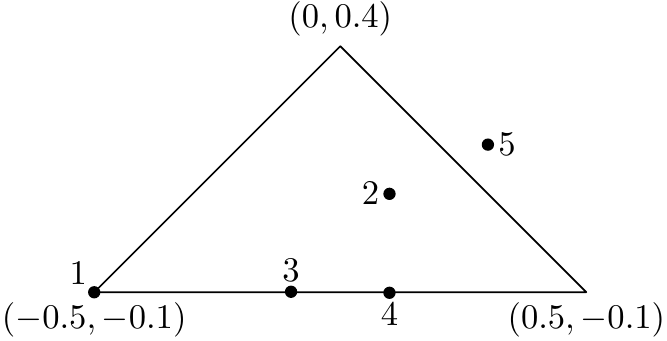}
\caption{Reference triangle for test of integration procedures. Points
  1--5 are the projections onto the triangle plane of the field
  points. Points are 1: $(-0.5,-0.1)$; 2: $(0.1,0.1)$; 3:
  $(0.1,-0.101)$; 4: $(-0.1,-0.099)$; 5: $(0.3,0.2)$.}
\label{fig:triangle}
\end{figure}

The triangular element for the first test is shown in
Figure~\ref{fig:triangle}. The test integrals were evaluated at a
range of values of $z$, over five different positions in the triangle
plane, shown in Figure~\ref{fig:triangle}. For $I_{1}$ and $I_{2}$,
$0\leq z\leq 1$, while for $I_{3}$, $1/8\leq z\leq 1$, since the
numerical integration routine will not give correct results at $z=0$.
These points were chosen to test for the cases most likely to be
encountered in practice. Point~1 is directly over a triangle vertex;
point~2 lies inside the triangle boundary; points~3 and~4 lie near the
triangle boundary but just inside and outside it respectively; point~5
is well separated from the element. Three basic integrals were
evaluated:
\begin{subequations}
  \label{equ:tests}
  \begin{align}
    I_{1} &= \iint \frac{L_{i}(\xi,\eta)}{4\pi R}\,\D x_{1}\,\D y_{1},\\
    I_{2} &= z\iint \frac{L_{i}(\xi,\eta)}{4\pi R^{3}}\,\D x_{1}\,\D
    y_{1} = -\frac{\partial I_{1}}{\partial n},\\
    I_{3} &= \iint \frac{L_{i}(\xi,\eta)}{4\pi R^{5}}\,\D x_{1}\,\D y_{1},
  \end{align}
\end{subequations}
where $L_{i}(\xi,\eta)$, $i=1,2,3$ are the linear shape functions for
the triangle:
\begin{align}
  \label{equ:shape}
  L_{1} &= 1-\xi-\eta,\quad L_{2} = \xi,\quad L_{3} = \eta.
\end{align}
For comparison with a numerical method, the integrals were also
evaluated using the approach of Hayami~\cite{hayami05}. The error
$\epsilon_{j}$, $j=1,2,3$, is calculated as the maximum of the
differences between the numerical and analytical values for $I_{j}$,
for each value of $i$, and is shown in Table~\ref{tab:triangle}. 

\begin{table}
  \caption{Errors in evaluation of integrals for triangle of
    Figure~\ref{fig:triangle}} 
  \label{tab:triangle}
  \centering
  
  \begin{tabular}{crrr}
    \hline
    Point & 
    \multicolumn{1}{c}{$\epsilon_{1}$} &
    \multicolumn{1}{c}{$\epsilon_{2}$} &
    \multicolumn{1}{c}{$\epsilon_{3}$} \\
    \hline
    1 & $1.5\times10^{-14}$ & $1.4\times10^{-12}$ & $4.4\times10^{-10}$ \\
    2 & $2.5\times10^{-11}$ & $1.4\times10^{-09}$ & $4.3\times10^{-07}$ \\
    3 & $4.9\times10^{-04}$ & $6.9\times10^{-09}$ & $2.4\times10^{-06}$ \\
    4 & $1.7\times10^{-05}$ & $6.8\times10^{-09}$ & $2.4\times10^{-06}$ \\
    5 & $4.6\times10^{-14}$ & $2.9\times10^{-12}$ & $8.9\times10^{-10}$ \\
    \hline
  \end{tabular}
\end{table}

From Table~\ref{tab:triangle}, the accuracy of the method is clearly
displayed: the errors are comparable to machine precision, with the
exception of the hypersingular case for points 2--4, where the
numerical integration routine would be expected to break down, and for
points~3 and~4 for $I_{1}$, where the numerical method has difficulty
in dealing with the very slender triangle between the evaluation point
and the lower vertices of the triangle. 

\subsection{Potential gradient calculation}
\label{sec:gradient}

\begin{figure}
  \centering
  \includegraphics[scale=0.5]{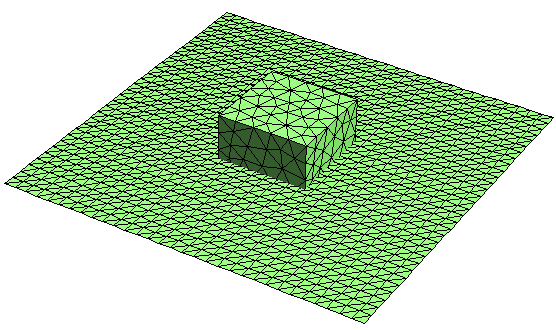}
  \caption{Unit cube and gradient test surface}
  \label{fig:cube}
\end{figure}

The second test was an evaluation of the gradient of a potential field
outside the unit cube. The evaluation method was implemented in the
free BEM code BEM3D~\cite{bem3d}, based on the GTS triangulated
surface library~\cite{popinet00}, and used to find the gradient of the
potential outside a unit cube centred on the origin,
Figure~\ref{fig:cube}. The boundary condition for potential and
potential gradient was imposed using a point source placed inside the
cube at $(1/4,1/4,1/4)$ and the field was evaluated as the gradient of
Equation~\ref{equ:potential}. The measure of error was the maximum
difference in any of the three components of the gradient on a
$32\times32$ grid $-2\leq x\leq2,-2\leq y\leq2,z=0$:
\begin{align}
  \label{equ:gradient:error}
  \epsilon &= 
  \max
  \left|
    \left.\frac{\partial\phi}{\partial x_{i}}\right|_{\text{computed}}
    -
    \left.\frac{\partial\phi}{\partial x_{i}}\right|_{\text{exact}}
  \right|
  /
  \max|\nabla\phi|.
\end{align}
The cube mesh was generated using
\textsc{gmsh}~\cite{geuzaine-remacle09}, with the discretization
length $h$ varied as a parameter.

\begin{figure}
  \centering
  \includegraphics{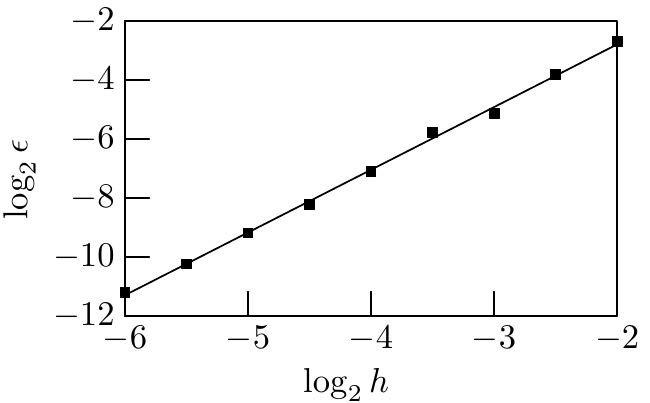}
  \caption{Error \textit{versus} discretization length for gradient of
    Laplace potential outside the unit cube. Symbols: computed error;
    solid line: $\epsilon=2.72h^{2.12}$.}
  \label{fig:error:cube}
\end{figure}

The resulting error is shown in Figure~\ref{fig:error:cube} and
demonstrates second order convergence with $h$, as might be expected
for linear elements. Again, the accuracy of the method is confirmed.

\section{Conclusions}
\label{sec:conclusions}

A method and explicit formulae for the exact integration of potential
integrals on planar triangular elements in the boundary element method
has been presented. The formulae are easily implemented using standard
operations and testing against numerical quadrature has shown them to
be accurate and, when implemented in a BEM code, convergent. 


\appendix

\section{Basic integrals}
\label{sec:basic}

Expressions for $I_{m,n}^{(p)}$ and $J_{m,n}$ can be found in standard
tables~\cite{gradshteyn-ryzhik80,prudnikov-brychkov-marichev03:1} and
are listed in Tables~\ref{tab:Imnp} and~\ref{tab:Jmn}. Also, where
$I_{m,n}^{(-3)}$ is not given explicitly, it can be expressed in terms
of lower order integrals~\cite[2.581.2]{gradshteyn-ryzhik80}:
\begin{align}
  \label{equ:Imn:3}
  I_{m,n}^{(-3)}(\alpha,\theta) &= 
    \frac{\sin^{m-1}\theta\cos^{n-1}\theta}{\alpha^{2}\Delta}
  -
  \frac{m-1}{\alpha^{2}}I_{m-2,n}^{(-1)}(\alpha,\theta)
  +
  \frac{n-1}{\alpha^{2}}I_{m,n-2}^{(-1)}(\alpha,\theta)
\end{align}

\begin{table}
  \caption{$I_{mn}^{(p)}(\theta)$, constant of integration omitted}
  \label{tab:Imnp}
  \centering
  
  \begin{tabular}{rrrr}
    \hline
    $p$ & $m$ & $n$ & 
    \multicolumn{1}{c}{$I_{mn}^{(p)}$}\\
    \hline
    -3 & 0 & 1 & 
    $\displaystyle \frac{\sin\theta}{\Delta}$
    \\[2ex]
    & 1 & 0 & 
    $\displaystyle -\frac{1}{\alpha'^{2}}\frac{\cos\theta}{\Delta}$
    \\[2ex]
    & 0 & 3 & 
    $\displaystyle
    -\frac{\alpha'^{2}}{\alpha^{2}}\frac{\sin\theta}{\Delta} + 
    \frac{1}{\alpha^{3}}\sin^{-1}(\alpha \sin\theta)$
    \\[2ex]
    -1 & 0 & 1 & 
    $\displaystyle \frac{1}{\alpha}\sin^{-1}(\alpha\sin\theta)$
    \\[2ex]
       & 1 & 0 & 
    $\displaystyle -\frac{1}{\alpha}\ln (\alpha \cos\theta + \Delta)$
    \\
       & 2 & -1 & 
    $\displaystyle \frac{1}{2\alpha'}
    \ln\frac{\Delta+\alpha'\sin\theta}{\Delta-\alpha'\sin\theta}
    - \frac{1}{\alpha}\sin^{-1}(\alpha\sin\theta)$
    \\[2ex]
    1  & 0 & -1 & 
    $\displaystyle \frac{\alpha'}{2}
    \ln\frac{\Delta+\alpha'\sin\theta}{\Delta-\alpha'\sin\theta}
    + \alpha\sin^{-1}(\alpha\sin\theta)$
    \\[2ex]
       & 2 & -1 & 
    $\displaystyle 
    -\frac{\Delta\sin\theta}{2} + 
    \frac{2\alpha^{2}-1}{2\alpha}\sin^{-1}(\alpha\sin\theta) +
    \frac{\alpha'}{2}
    \ln\frac{\Delta+\alpha'\sin\theta}{\Delta-\alpha'\sin\theta}
    $
    \\[2ex]
       & 1 & -2 & 
    $\displaystyle \frac{\Delta}{\cos\theta} - 
    \alpha\ln(\alpha\cos\theta + \Delta)$
    \\[2ex]
       & 1 & 0 & 
    $\displaystyle 
    -\frac{\Delta\cos\theta}{2} - 
    \frac{1}{2}\frac{\alpha'^{2}}{\alpha}
    \ln(\alpha\cos\theta + \Delta)
    $
    \\
    \hline
  \end{tabular}
\end{table}

\begin{table}
  \caption{$J_{mn}(\theta)$, constant of integration omitted}
  \label{tab:Jmn}
  \centering
  
  \begin{tabular}{rrr}
    \hline
    $m$ & $n$ & \multicolumn{1}{c}{$J_{mn}$} \\
    \hline
    0 & 0 & $\theta$ \\[2ex]
    0 & -1 &
    $\displaystyle 
    \ln
      \frac{1+\sin\theta}{1-\sin\theta}
    $
    \\[2ex]
    0 & 1 &
    $\displaystyle \sin\theta$
    \\[2ex]
    1 & 0 &
    $\displaystyle -\cos\theta$
    \\[2ex]
    0 & 3 &
    $\displaystyle \sin\theta - \frac{1}{3}\sin^{3}\theta$
    \\[2ex]
    1 & -2 &
    $\displaystyle \frac{1}{\cos\theta}$
    \\[2ex]
    1 & 2 &
    $\displaystyle -\frac{\cos^{3}\theta}{3}$
    \\[2ex]
    2 & -1 &
    $\displaystyle -\sin\theta + 
    \ln
    \left(
      \frac{\pi}{4} + \frac{\theta}{2}
    \right)$
    \\[2ex]
    \hline
  \end{tabular}
\end{table}


\end{document}